\documentclass[12pt]{article}
\usepackage{amsmath,amsthm,amssymb}
\usepackage{mathtools}

\newtheorem{theorem}{Theorem}[section]

\newtheorem{remark}{Remark}[section]

\numberwithin{equation}{section}

\newcommand{\D}{{\rm d}}
\newcommand{\dx}{\, \D x}

\newcommand{\dt}{\, \D t}

\newcommand{\rz}{\mathbb{R}}

\newcommand{\eps}{\varepsilon}
\newcommand{\iom}{\int_{\Omega}}

\newcommand{\psp}{\,}

\title{Variants of Bernstein's theorem for variational integrals with linear and nearly linear growth}

\author{Michael Bildhauer \& Martin Fuchs}
\date{}

\usepackage{hyperref}

\hypersetup{
   pdfnewwindow=true,
   colorlinks=false,
   linkcolor=black,
   citecolor=green,
   filecolor=magenta,
   urlcolor=cyan,
}

\begin{document}

\parindent0em
\maketitle

\newcommand{\op}[1]{\operatorname{#1}}
\newcommand{\bv}{\op{BV}}
\newcommand{\mub}{\overline{\mu}}
\newcommand{\muhat}{\hat{\mu}}

\newcommand{\hypref}[2]{\hyperref[#2]{#1 \ref*{#2}}}
\newcommand{\hypreff}[1]{\hyperref[#1]{(\ref*{#1})}}

\newcommand{\ob}[1]{^{(#1)}}

\newcommand{\xh}{\Xi}
\newcommand{\oh}[1]{O\left(#1\right)}
\newcommand{\xn}{w_{1}}
\newcommand{\yn}{w_{2}}
\newcommand{\On}{\hat{\Omega}}
\newcommand{\uq}{\hat{u}}
\newcommand{\gameps}{\Gamma_{Q}}
\newcommand{\se}{{s_Q}}
\begin{abstract}{\footnote{AMS subject classification: 49Q20, 49Q05, 53A10, 35J20\\
Keywords: Bernstein's theorem, non-parametric minimal surfaces, variational problems with (nearly) linear growth, equations in two variables}}
Using a Caccioppoli-type inequality involving negative exponents for a directional weight  we establish variants of Bernstein's theorem
for variational integrals with linear and nearly linear growth. We give some mild conditions for entire solutions of the equation
\[
\op{div} \Big[Df(\nabla u)\Big] = 0 \, ,
\]
under which solutions have to be affine functions. Here $f$ is a smooth energy density satisfying $D^2 f>0$
together with a natural growth condition for $D^2 f$.
\end{abstract}

\newcommand{\platz}{\vspace{2ex}}
\newcommand{\absatz}{
\platz
\centerline{\rule[1ex]{3cm}{0.01cm}}
\platz
}

\newcommand{\os}{\overline{s}_\gamma}
\newcommand{\oslimit}{\overline{s}}
\newcommand{\ls}{\varkappa}

\newcommand{\gamf}{{\theta_{i}}}
\newcommand{\gamg}{{\gamma_{i}}}

\parindent0ex
\section{Introduction}\label{in}


In \cite{Be:1927_1} Bernstein proved that every $C^2$-solution $u=u(x) = u(x_1,x_2)$ of the non-parametric
minimal surface equation
\begin{equation}\label{in 1}
\op{div} \Bigg[\frac{\nabla u}{\sqrt{1+|\nabla u|^2}}\Bigg] = 0
\end{equation}
over the entire plane must be an affine function, which means that with real numbers $a$, $b$, $c$ it holds
\[
u(x_1,x_2) = a x_1 + b x_2 + c\, .
\]
For a detailed discussion of this classical result the interested reader is
referred for instance to \cite{DHS:2010_1}, \cite{Ni:1989_1}, \cite{Ni:1957_1}, \cite{Os:1986_1} and the references quoted therein.\\

Starting from Bernstein's result the question arises to which classes of second order equations the Bernstein-property
extends. More precisely, we replace \eqref{in 1} through the equation
\begin{equation}\label{in 2}
L[u] = 0
\end{equation}
for a second order elliptic operator $L$ and assume that $u \in C^2(\mathbb{R}^2)$ is an entire solution of \eqref{in 2}
asking if $u$ is affine. To our knowledge a complete answer to this problem is open, however we have the explicit
``Nitsche-criterion'' established by J.C.C.~Nitsche and J.A.~Nitsche \cite{NN:1959_1}.\\

In our note we discuss equation \eqref{in 2} assuming that $L$ is the Euler-operator associated to the variational integral
\[
J[u,\Omega] := \iom f(\nabla u) \dx
\]
with density $f$: $\rz^2 \to \rz$ and for domains $\Omega \subset \rz^2$, i.e. \eqref{in 2} is replaced by
\begin{equation}\label{in 3}
\op{div} \Big[Df(\nabla u)\Big] = 0 \, ,
\end{equation}
where in the minimal surface case \eqref{in 1} we have $f(p) = \sqrt{1+|p|^2}$, $p\in \rz^2$, being an integrand of linear growth with repect 
to $\nabla u$. 
For this particular class of energy densities and under the additional assumption that $f$ is of type $f(p) = g(|p|)$, $p\in \rz^2$, for
a function $g\in C^2([0,\infty))$ such that
\[
0 < g''(t) \leq c (1+t)^{-\mu}\, , \quad t \geq 0\, ,
\]
with exponent $\mu \geq 3$ (including the minimal surface case) we proved the Bernstein-property in Theorem 1.2 of \cite{BF:2021_1} benefiting
from the work \cite{FSV:2008_1} of Farina, Sciunzi and Valdinoci.\\ 

Bernstein-type theorems under natural additional conditions to be imposed on the entire solutions of the Euler-equations for splitting-type variational integrals of linear growth
have been established in the recent paper \cite{BFF:2023_1}.\\

One of the main tools used in \cite{BFF:2023_1} is a Caccioppoli-inequality involving negative exponents which was already exploited in 
different variants in the papers \cite{BF:2020_3}, \cite{BF:2020_4}, \cite{BF:2021_3}, \cite{BF:2022_1}.\\

In fact we used this inequality to show that $\partial_1 \nabla u \equiv 0$ which follows by considering the bilinear form $D^2f(\nabla u)$  with suitable weights.
Since we are in two dimensions, the use of equation \eqref{in 3} then completes the proof of the splitting-type results in \cite{BFF:2023_1}.\\

In the manuscript at hand we observe that even without splitting-structure it is possible to discuss the Caccioppoli-inequality with a directional weight
obtaining $\partial_1 \nabla u \equiv 0$ and to argue similar as before. We note that considering directional weights gives much more flexibility in choosing 
exponents than arguing with a full gradient (see Proposition 6.1 and Proposition 6.2 of \cite{BF:2022_1}). We here already note that we also include
a logarithmic variant of the Caccioppli-inequality as an approach to the limit case $\alpha = -1/2$.\\

Before going into these details let us have a brief look at power growth energy densities as for example
\[
f(p) = (1+|p|^2)^{\frac{s}{2}} \, , \quad p \in \rz^2\, ,
\]
with exponent $s >1$. Then the Nitsche-criterion (compare \cite{NN:1959_1}, Satz) shows the existence of non-affine
entire solutions to equation \eqref{in 3}, and as we will shortly discuss in the Appendix the same reasoning applies
to the nearly linear growth model
\begin{equation}\label{in 4}
f(p) = |p| \ln (1+|p|)\, , \quad p \in \rz^2\, ,
\end{equation}
which means that we do not have the Bernstein-property for equation \eqref{in 3} with density of the form \eqref{in 4}. \\

However, as indicated above, we can establish a mild condition under which any entire solution is an affine function being valid for a
large class of densities $f$ including the nearly linear and even the linear case.\\

Let us formulate our\\ 

{\bf Assumptions.} The density $f$: $\rz^2 \to \rz$ is of class $C^2$ such that 
\begin{equation}\label{in 5}
D^2f(p)(q,q) >0 \, \quad \mbox{for all $p$, $q\in \rz^2$, $q\not= 0$.}
\end{equation}

For a constant $\lambda >0$ it holds
\begin{equation}\label{in 6}
D^2 f(p)(q,q)\leq \lambda \frac{\ln(2+|p|}{1+|p|}|q|^2 \, , \quad \mbox{$p$, $q \in \rz^2$}.
\end{equation}
\begin{remark}\label{rem 1 1}
Condition \eqref{in 5} implies the strict convexity of $f$, whereas from \eqref{in 6} we obtain
\[
|f(p)| \leq c \big( |p| \ln (1+|p|) + 1\big)\, , \quad p \in \rz^2\, , 
\]
with some constant $c>0$.
\end{remark}

We have the following result:
\begin{theorem}\label{in theo 1}
Let $f$ satisfy \eqref{in 5} and \eqref{in 6} and consider an entire solution $u \in C^2(\rz^2)$
of equation \eqref{in 3}. Suppose that  with numbers $0\leq m < 1$, $K >0$ the solution satisfies 
\begin{equation}\label{in 7}
|\partial_1 u(x)| \leq K \Big( |\partial_2 u(x)|^m + 1 \big)\, , \quad x \in \rz^2\, ,
\end{equation}
or
\begin{equation}\label{in 8}
|\partial_2 u(x)| \leq K \Big( |\partial_1 u(x)|^m + 1 \big)\, , \quad x \in \rz^2\, .
\end{equation}
Then $u$ is affine.
\end{theorem}

\begin{remark}\label{rem 1 2}
Since the density $f$ from \eqref{in 4} fulfills the the conditions \eqref{in 5} and \eqref{in 6} and since in this case non-affine
entire solutions exist, the requirements \eqref{in 7} and \eqref{in 8} single out a class of entire solutions of Bernstein-type.

Of course we know nothing concerning the optimality of \eqref{in 7} and \eqref{in 8}. Another unsolved problem is the question,
if in the case of linear growth with radial structure, i.e.~$f(\nabla u) = g(|\nabla  u|)$, Bernstein's theorem holds without extra conditions
on the entire solution $u$.
\end{remark}

\begin{remark}\label{rem 1 3}
The conditions \eqref{in 7} and \eqref{in 8} are in some sense related to the  ``balancing conditions'' used in \cite{BFF:2023_1}
in order to exclude entire solutions of the form $u(x_1,x_2) = x_1x_2$ for densities $f$ of splitting type. 
\end{remark}

Let us pass to the linear growth case replacing \eqref{in 6} by
\begin{equation}\label{in 9}
|D^2f(p)| \leq \Lambda \frac{1}{1+|p|}\, , \quad p\in \rz^2\, ,
\end{equation}
with a positive constant $\Lambda$. Here the notion of linear growth just expresses the fact that from \eqref{in 9} it follows
that
\[
|f(p)| \leq c \big(|p| +1\big) \, , \quad p \in \rz^2\, ,
\]
with some number $c >0$. In this situation we have
\begin{theorem}\label{in theo 2}
Let $f$ satisfy \eqref{in 5} together with \eqref{in 9} and let $u \in C^2(\rz^2)$ denote an entire solution of equation 
\eqref{in 3} for which we have
\begin{equation}\label{in 10}
|\partial_1 u(x)| \ln^2\big(1+|\partial_1 u(x)|\big) \leq K \big(|\partial_2 u(x)| +1\big) \, ,\quad x \in \rz^2\, ,
\end{equation}
or
\begin{equation}\label{in 11}
|\partial_2 u(x)| \ln^2\big(1+|\partial_2 u(x)|\big) \leq K \big(|\partial_1 u(x)| +1\big) \, ,\quad x \in \rz^2\, ,
\end{equation}
with some number $K \in (0,\infty)$. Then $u$ is an affine function.
\end{theorem}

The results of Theorem \ref{in theo 1} and Theorem \ref{in theo 2} are not limited to the particular coordinate directions
$e_1=(1,0)$ and $e_2 =(0,1)$, more precisely it holds:

\begin{theorem}\label{in theo 3}
Let $f$ satisfy either the assumptions of Theorem \ref{in theo 1} (``case 1'') or of Theorem \ref{in theo 2} (``case 2'')
and suppose that $u \in C^2(\rz^2)$ is an entire solution of \eqref{in 3}. Assume that there exist two linearly independent
vectors $E_1$, $E_2 \in \rz^2$ such that\\

\begin{tabular}{ll}
in case 1:& \eqref{in 7} holds with $\partial_\alpha u$ being replaced by $\partial_{E_\alpha} u$, $ \alpha=1$, $2$,\\[2ex]
in case 2:& \eqref{in 10} is true again with $\partial_\alpha u$ being replaced by $\partial_{E_\alpha} u$, $\alpha =1$, $2$.
\end{tabular} 
Then $u$ is an affine function.
\end{theorem} 
The proof of this result follows from the observation that the function $u$ is a local minimizer of the energy
$\int f(\nabla u)\dx$ combined with a suitable linear transformation. If we let
\begin{eqnarray*}
T:&& \rz^2\to \rz^2\, , \quad E_\alpha = T(e_\alpha)\, , \quad \alpha = 1, \, 2\, ,\\[2ex]
E_{\alpha\beta}&:=&E_\alpha \cdot E_\beta\, , \quad \alpha ,\, \beta = 1,\, 2\, ,\\[2ex]
\tilde{u}(x) &:=& u\big(T(x)\big)\, ,\quad x\in \rz^2 \, ,\\[2ex]
\tilde{f}(p)&:=& f\Bigg(T\Big(\big(E_{\alpha \beta}\big)^{-1}_{1\leq \alpha , \beta \leq 2}\, p\Big)\Bigg)\, , \quad p \in \rz^2\, ,
\end{eqnarray*}
then it holds
\[
\partial_\alpha \tilde{u}(x) = \partial_{E_\alpha} u \big(T(x)\big)\, , \quad x \in \rz^2\, , \quad \alpha =1,\, 2\, ,
\]
and $\tilde{u}$ is an entire solution of equation \eqref{in 3} with $f$ being replaced by $\tilde{f}$, which follows from the local minimality
of $\tilde{u}$ with respect to the energy $\int \tilde{f}(\nabla w) \dx$. Obviously the properties of $\tilde{f}$ required in Theorem \ref{in theo 1}
and Theorem \ref{in theo 2}, respectively, are consequences of the corresponding assumptions imposed on $f$, thus we can apply our previous results
to $\tilde{u}$ (and $\tilde{f}$).\\

Our paper is organized as follows: in Section \ref{pr 1} we present the proof of Theorem \ref{in theo 1} based on a Caccioppoli-inequality
involving negative exponents, which has been established, for instance, in \cite{BF:2022_1}, Proposition 6.1.\\

Section \ref{pr 2} is devoted to the discussion of Theorem \ref{in theo 2}. We will make use of some kind of a limit version of
Caccioppoli's inequality, whose proof will be presented below.
With the help of this result the claim of Theorem \ref{in theo 2} follows along the lines of Section 2.
We finish Section 3 by presenting a technical extension of Theorem \ref{in theo 2}, which just follows from an inspection
of the arguments (compare Theorem \ref{pr 2 theo 1}).\\

For the reader's convenience we discuss in an appendix equation \eqref{in 3} for the nearly linear growth case \eqref{in 4} and
show that the Nitsche-criterion applies yielding non-affine solutions defined on the whole plane.

\section{Proof of Theorem \ref{in theo 1}}\label{pr 1}

Let $f$ satisfy \eqref{in 5} and \eqref{in 6}, let $u$ denote an entire solution of \eqref{in 3} and assume w.l.o.g.~that
\eqref{in 7} holds. We apply inequality (107) from Proposition 6.1 in \cite{BF:2022_1} with the choices $l=1$, $i=1$ and 
\[
\Omega = B_{2R}= \big\{x\in \rz^2: \, |x| < 2R\big\}
\]
to obtain for any $\alpha > -1/2$ and $\eta \in C^\infty_0(B_{2R})$, $0\leq \eta < 1$,
\begin{eqnarray}\label{pr 1 1}
\lefteqn{ \int_{B_{2R}} \eta^2 D^2f(\nabla u)\big(\nabla \partial_1 u, \nabla \partial_1 u\big) \Gamma_1^\alpha \dx}\nonumber\\
&\leq & 
c \int_{B_{2R}} D^2f(\nabla u) (\nabla \eta,\nabla \eta) \Gamma_1^{\alpha+1} \dx\, , \quad \Gamma_1 := 1+ |\partial_1 u|^2\, ,  
\end{eqnarray}
with a finite constant independent of $R$.\\

Letting $\eta=1$ on $B_R$ and assuming $|\nabla \eta| \leq c/R$ we apply \eqref{in 6} to the r.h.s.~of \eqref{pr 1 1} and get
\begin{eqnarray}\label{pr 1 2}
\lefteqn{\int_{B_R} D^2f(\nabla u)\big( \nabla \partial_1 u,\nabla \partial_1 u) \Gamma_1^\alpha \dx}\nonumber\\
&\leq & cR^{-2} \int_{R < |x|<2R} \Gamma_1^{1+\alpha} \ln\big(2+|\nabla u|\big)\big(1+|\nabla u|\big)^{-1} \dx \, .
\end{eqnarray}
On account of \eqref{in 7} we deduce for any $\eps >0$
\begin{eqnarray*}
\Gamma_1^{1+\alpha} \ln\big(2+|\nabla u|\big)\big(1+|\nabla u|\big)^{-1} & \leq &
c(\eps) \big( 1+|\partial_2 u|^2\big)^{m(1+\alpha)} \big(1+|\nabla u|\big)^{\eps - 1}\\[2ex]
&\leq & \tilde{c}(\eps) \big(1+|\partial_ 2 u|\big)^{2m(1+\alpha)} \big(1+|\partial_2 u|\big)^{\eps-1}\\[2ex]
&=& \tilde{c}(\eps) \big(1 + |\partial_2 u|\big)^{2m(1+\alpha)-1+\eps} \, .
\end{eqnarray*}

Recall that $m < 1$, hence $2m(1+\alpha) -1 <0$ for $\alpha > -1/2$ sufficiently close to $-1/2$. We fix $\alpha$
with this property and finally select $\eps >0$ such that $2m(1+\alpha) - 1+ \eps \leq 0$ to obtain
\begin{equation}\label{pr 1 3}
\Gamma_1^{1+\alpha} \ln\big(2+|\nabla u|\big) \big(1+|\nabla u|\big)^{-1} \leq const < \infty \, .
\end{equation}
Combining \eqref{pr 1 2} with \eqref{pr 1 3} it is shown that
\begin{equation}\label{pr 1 4}
\int_{\rz^2} D^2f(\nabla u) \big(\nabla \partial_1 u,\nabla \partial_1 u\big) \Gamma_1^\alpha \dx < \infty \, .
\end{equation}
We quote equation (108) from \cite{BF:2022_1} again with the previous choices $l=1$, $i=1$, $\Omega = B_{2R}$,
$\eta \in C^\infty_0(B_{2R})$, $0\leq \eta \leq 1$, and with $\alpha$ as fixed above. The same calculations
as carried out after (108) then yield
\begin{eqnarray}\label{pr 1 5}
\lefteqn{\int_{B_{2R}} D^2f(\nabla u) \big(\nabla \partial_1 u,\nabla \partial_1 u\big) \eta^2 \Gamma_1^\alpha \dx}\nonumber\\[2ex]
&\leq & c \Bigg| \int_{B_{2R}-B_R} D^2f(\nabla u) \big(\nabla \partial_1 u,\nabla \eta^2\big) \partial_1 u \Gamma_1^\alpha \dx \Bigg|\, .
\end{eqnarray}
On the r.h.s.~of \eqref{pr 1 5} we apply the Cauchy-Schwarz inequality to the bilinear form $D^2f(\nabla u)$, hence
\begin{eqnarray}\label{pr 1 6}
\mbox{l.h.s.~of \eqref{pr 1 5}}& \leq &
\Bigg[\int_{B_{2R}-B_R} D^2f(\nabla u) \big(\nabla \partial_1 u, \nabla \partial_1 u\big) \Gamma_1^\alpha \eta^2 \dx\Bigg]^{\frac{1}{2}}
\nonumber\\[2ex]
&& \quad \cdot \Bigg[\int_{B_{2R}-B_R} D^2f(\nabla u) (\nabla \eta,\nabla \eta) \Gamma_1^{1+\alpha}\dx\Bigg]^{\frac{1}{2}}\nonumber\\[2ex]
& =: & T_1(R)^{\frac{1}{2}} \cdot T_2(R)^{\frac{1}{2}} \psp .
\end{eqnarray}
Here $\eta$ has been chosen in such a way that $\eta \equiv 1$ on $B_R$ and therefore 
$\op{spt}(\nabla \eta) \subset B_{2R} -B_R$. By \eqref{pr 1 4} we have
\[
\lim_{R \to \infty} T_1(R) = 0 \, ,
\]
whereas the calculations carried out after \eqref{pr 1 2} imply the boundedness of $T_2(R)$. Thus \eqref{pr 1 5} and \eqref{pr 1 6}
imply
\[
\int_{\rz^2} D^2f(\nabla u) \big( \nabla \partial_1 u , \nabla \partial_1 u\big) \Gamma_1^\alpha \dx = 0 \, ,
\]
hence $\nabla \partial_1 u = 0$ on account of \eqref{in 5}. This shows $\partial_1 u =a$ for some number $a \in \rz$
and since
\[
u(x_1,x_2) - u(0,x_2) = \int_0^{x_1} \frac{\D}{\dt} u(t,x_2) \dt = a x_1
\]
we can write 
\[
u(x_1,x_2) = \varphi(x_2) + a x_1\, ,\quad \varphi(x_2) := u(0,x_2) \, .
\]
Equation \eqref{in 3} gives
\[
\frac{\D}{\D t} \frac{\partial f}{\partial p_2}\big(a,\varphi'(t)\big) = 0 
\]
so that
\[
\frac{\partial f}{\partial p_2} \big(a, \varphi'(t)\big) = c
\]
for a constant $c$. Finally we observe that the function
\[
y \mapsto \frac{\partial f}{\partial p_2} (a,y)
\]
is strictly increasing (recall \eqref{in 5}), which shows the constancy of $\varphi'$ and therefore
\[
u(0,x_2) = b x_2 +c
\]
for some numbers $b$, $c \in \rz$. Altogether we have shown that $u$ is affine finishing the proof of
Theorem \ref{in theo 1}. \hspace*{\fill}\qed\\

\section{Proof of Theorem \ref{in theo 2}}\label{pr 2}

Let the assumptions of Theorem \ref{in theo 2} hold and consider an entire solution $u \in C^2(\rz^2)$ of \eqref{in 3} without
requiring \eqref{in 10} or \eqref{in 11} for the moment. If we use condition \eqref{in 9} in inequality \eqref{pr 1 1} and if we assume 
that the choice $\alpha = -1/2$ is admissible in \eqref{pr 1 1}, then the calculations of Section \ref{pr 1} would immediately imply 
that $\nabla^2 u =0$ yielding Bernstein's theorem, i.e.~the entire solution $u$ is an affine function without adding further hypotheses on $u$.\\

However, we do not have \eqref{pr 1 1} in the case that $\alpha = -1/2$ and hence we provide a weaker version involving conditions like
\eqref{in 10} or \eqref{in 11} in order to conclude that $u$ is affine.\\

To be precise, we assume the validity of \eqref{in 10}, let $l$, $i$, $\Omega$ and $\eta$ as stated in front of \eqref{pr 1 1} recalling
\[
\int_{B_{2R}} D^2f(\nabla u) \big(\nabla \partial_1 u, \nabla \psi\big) \dx = 0
\]
for the choice $\psi: = \eta^2 \partial_1 u \Phi(\Gamma_1)$, where
\begin{equation}\label{pr 2 1}
\Phi(t) := \frac{\ln(e^2-1+t)}{\sqrt{t}} \, , \quad t \geq 1\, .
\end{equation}
Note that the choice $\alpha = -1/2$ is compensated by the logarithm. Obviously $\Phi(1) =2$, $\Phi(\infty) = 0$ together with
\begin{equation}\label{pr 2 2}
\Phi'(t) = -\frac{1}{2} \frac{1}{t^{\frac{3}{2}}} \ln(e^2-1+t) + \frac{1}{\sqrt{t} (e^2-1+t)} < 0\, , \quad t \geq 1\, ,
\end{equation}
where the negative sign for $\Phi'(t)$ follows from $\ln(e^2 -1 +t) \geq 2$ for $t \geq 1$.\\

With $\psi$ from above and $\Phi$ defined according to \eqref{pr 2 1} we obtain
\begin{eqnarray}\label{pr 2 3}
\lefteqn{ \int_{B_{2R}} \eta^2 D^2f(\nabla u)\big(\nabla \partial_1 u,\nabla \partial_1 u\big) \Phi(\Gamma_1) \dx}\nonumber\\[2ex]
&+& \int_{B_{2R}} \eta^2 D^2 f(\nabla u)\big(\nabla \partial_1 u,\partial_1 u \nabla \Phi(\Gamma_1)\big)\dx\nonumber\\[2ex]
&&= -2 \int_{B_{2R}} \eta D^2f(\nabla u) \big(\nabla \partial_1 u, \nabla \eta\big) \partial_1 u \Phi(\Gamma_1) \dx\, .
\end{eqnarray}
Using the identity
\begin{eqnarray*}
\partial_1 u \nabla \Phi(\Gamma_1) &=& 2 \nabla \partial_1 u (\partial_1 u)^2 \Phi'(\Gamma_1)\\[2ex]
&=& 2 \nabla \partial_1 u \Bigg[ \Gamma_1 \Phi'(\Gamma_1) - \Phi'(\Gamma_1)\Bigg]\, ,
\end{eqnarray*}
the left-hand side of \eqref{pr 2 3} equals
\[
\int_{B_{2R}} \eta^2 D^2f(\nabla u) \big(\nabla \partial_1 u,\nabla \partial_1 u\big) \Big[ \Phi(\Gamma_1) 
+ 2 \Gamma_1 \Phi'(\Gamma_1) - 2 \Phi'(\Gamma_1)\Big]\dx \, .
\]
From \eqref{pr 2 2} it follows (recalling $\Phi'(t) \leq 0$ for $t \geq 1$)
\begin{eqnarray*}
\lefteqn{\Phi(\Gamma_1)  + 2 \Gamma_1 \Phi'(\Gamma_1) - 2 \Phi'(\Gamma_1)}\\[2ex]
&\geq &
\Phi(\Gamma_1) + 2 \Gamma_1 \Phi'(\Gamma_1)\\[2ex]
&=& \frac{\ln(e^2-1+\Gamma_1)}{\sqrt{\Gamma_1}}
+ 2 \Gamma_1 \Bigg[-\frac{1}{2}\frac{\ln(e^2-1+\Gamma_1)}{\Gamma_1^{\frac{3}{2}}} + \frac{1}{\sqrt{\Gamma_1} (e^2-1+\Gamma_1)}\Bigg]\\[2ex]
&=& 2 \frac{\sqrt{\Gamma_1}}{e^2-1+\Gamma_1}  \geq c \frac{1}{\sqrt{\Gamma_1 }} 
\end{eqnarray*}
for some constant $c >0$. Altogether we deduce from \eqref{pr 2 3} the inequality of Caccioppoli-type
\begin{eqnarray}\label{pr 2 4}
\lefteqn{\int_{B_{2R}} \eta^2 D^2f(\nabla u) \big(\nabla \partial_1 u, \nabla \partial_1 u\big) \Gamma_1^{-\frac{1}{2}} \dx}\nonumber\\[2ex]
&\leq & - 2 \int_{B_{2R}} \eta D^2f(\nabla u) \big(\nabla \partial_1 u, \nabla \eta\big) \partial_1 u \Gamma_1^{-\frac{1}{2}}\ln(e^2-1+\Gamma_1) \dx\nonumber\\[2ex] 
&=:& -2 S\, .
\end{eqnarray}
To the quantity $S$ we apply the Cauchy-Schwarz inequality valid for the bilinear form $D^2f(\nabla u)$ and get
\begin{eqnarray}\label{pr 2 5}
|S| & \leq & \int_{B_{2R}} \Bigg|D^2f(\nabla u) \big(\eta \Gamma_1^{-\frac{1}{4}} \nabla \partial_1 u,\partial_1 u \Gamma_1^{-\frac{1}{4}} \ln(e^2-1+\Gamma_1)\nabla \eta\big)\Bigg|\dx
\nonumber\\[2ex]
&\leq &\Bigg[\int_{B_{2R}} D^2f(\nabla u) \big(\nabla \partial_1 u,\nabla \partial_1 u\big)\eta^2 \Gamma_1^{-\frac{1}{2}} \dx\Bigg]^{\frac{1}{2}}\nonumber\\[2ex]
&&\qquad \cdot \Bigg[\int_{B_{2R}} D^2f(\nabla u)(\nabla \eta,\nabla \eta) \Gamma_1^{\frac{1}{2}} \ln^2(e^2-1+\Gamma_1)\dx\Bigg]^{\frac{1}{2}} \, .
\end{eqnarray} 
On the right-hand side of \eqref{pr 2 5} we make use of Young's inequality yielding for any $\eps >0$
\begin{eqnarray}\label{pr 2 6}
|S| & \leq & \eps \int_{B_{2R}} \eta^2 \Gamma_1^{-\frac{1}{2}} D^2f(\nabla u)\big(\nabla \partial_1 u,\nabla \partial_1 u\big) \dx\nonumber\\[2ex]
&&+ c(\eps) \int_{B_{2R}} D^2f(\nabla u) (\nabla \eta,\nabla \eta) \Gamma_1^{\frac{1}{2}} \ln^2(e^2-1+\Gamma_1) \dx \, .
\end{eqnarray}
Finally we combine \eqref{pr 2 6} and \eqref{pr 2 4}, thus for a fixed $\eps$ being sufficiently small it holds
\begin{eqnarray}\label{pr 2 7}
\lefteqn{\int_{B_{2R}} D^2f(\nabla u)\big(\nabla \partial_1 u,\nabla \partial_1 u\big)\eta^2 \Gamma_1^{-\frac{1}{2}}\dx}\nonumber\\[2ex]
&\leq & c \int_{B_{2R}} D^2f(\nabla u) (\nabla \eta,\nabla \eta) \Gamma_1^{\frac{1}{2}} \ln^2(e^2-1+\Gamma_1)\dx \, .
\end{eqnarray}
The properties of $\eta$ as stated after \eqref{pr 1 1} imply that the right-hand side of \eqref{pr 2 7} is bounded by (recall \eqref{in 9})
\begin{eqnarray*}
\lefteqn{cR^{-2} \int_{B_{2R}} |D^2f(\nabla u)| \Gamma_1^{\frac{1}{2}} \ln^2(e^2-1+\Gamma_1) \dx}\\[2ex]
&\leq & c R^{-2} \int_{B_{2R}} \frac{1}{1+|\nabla u|} \Gamma_1^{\frac{1}{2}} \ln^2(e^2-1+\Gamma_1) \dx\nonumber\\[2ex]
&\leq & c R^{-2} \int_{B_{2R}} \frac{1}{1+|\partial_2 u|} (1+|\partial_1 u|) \ln^2(1+|\partial_1 u|) \dx \, .
\end{eqnarray*}
Quoting \eqref{in 10} and returning to \eqref{pr 2 7} we find that
\begin{equation}\label{pr 2 8}
\int_{\rz^2} D^2f(\nabla u) \big(\nabla \partial_1 u,\nabla \partial_1 u\big) \eta^2 \Gamma_1^{-\frac{1}{2}} \dx < \infty \, .
\end{equation}
With \eqref{pr 2 8} and on account of estimate \eqref{pr 2 5} it is immediate (recall \eqref{pr 2 4}) that actually
\[
\int_{\rz^2} D^2f(\nabla u) \big(\nabla \partial_1 u,\nabla \partial_1 u\big) \Gamma_1^{-\frac{1}{2}} = 0 \, ,
\]
hence $\nabla \partial_1 u =0$ and we can follow the lines of Section \ref{pr 1} to prove our claim. \qed\\

An inspection of our previous arguments shows that we can replace the function $\ln(e^2-1+t)$ used before by any function
$\rho$: $[1,\infty) \to \rz^+$ of class $C^1$ such that we have for all $t\geq 1$
\begin{equation}\label{pr 2 9}
\rho'(t) >0\, ,\quad \frac{\D}{\D t} \Big[\frac{1}{\sqrt{t}} \rho(t) \Big] \leq 0 \, .
\end{equation}
Replacing \eqref{pr 2 1} by
\begin{equation}\label{pr 2 10}
\Phi(t) := \frac{1}{\sqrt{t}} \rho(t)\, \quad t\geq 1\, ,
\end{equation}
and letting as before $\psi := \eta^2 \partial_1 u \Phi(\Gamma_1)$ now with $\Phi$ defined in \eqref{pr 2 10}
we obtain
\begin{theorem}\label{pr 2 theo 1}
Let $f$ satisfy \eqref{in 5} together with \eqref{in 9} and choose $\rho$ according to \eqref{pr 2 9}. Suppose that $u \in C^2(\rz^2)$
is an entire solution of \eqref{in 3} such that
\[
\Gamma_1^{-\frac{1}{2}}\,  \frac{\rho^2(\Gamma_1)}{\rho'(\Gamma_1)}\leq c \Gamma_2^{\frac{1}{2}}\, , \quad
\Gamma_i := 1+|\partial_i u|^2\, , \quad i=1,\, 2\, ,
\]
or 
\[
\sup_{R>0} R^{-2} \int_{B_R} \Gamma_1^{-1} \, \frac{\rho^2(\Gamma_1)}{\rho'(\Gamma_1)} \dx < \infty
\]
holds with some finite constant $c$. Then $u$ is an affine function.
\end{theorem}
We leave the details to the reader just adding the obvious remark that clearly we can interchange the roles of the partial 
derivatives $\partial_1 u$ and $\partial_2 u$ or even work with arbitrary directional derivatives as done in Theorem \ref{in theo 3}.

\section{Appendix}\label{app}

We shortly discuss the Nitsche-criterion (see \cite{NN:1959_1}, Satz) for the model case
\[
\int |\nabla u|\ln\big(1+|\nabla u|\big) \dx \, .
\]
With a slight abuse of notation but in accordance with the terminology of \cite{NN:1959_1} we let
\[
g(t) := t \ln(1+t)\, , \quad t \geq 0\, ,\qquad f(t) := g\big(\sqrt{t}\big) \, ,
\]
so that
\[
J[u] := \int |\nabla u| \ln\big(1+|\nabla u|\big) \dx = \int f\big(|\nabla u|^2\big)\dx \, .
\]
Introducing the function $\lambda(t) := 2 f''(t)/f'(t)$ again for $t \geq 0$ we claim
\begin{equation}\label{app 1}
\int_1^\infty \frac{1+t \lambda(t)}{2+t \lambda(t)} \frac{1}{t} \dt = \infty \, .
\end{equation} 
From \eqref{app 1} it follows that the Euler-equation associated to the functional $J$ admits entire non-affine solutions.\\

For \eqref{app 1} we observe the formula
\[
\frac{1}{t} \frac{1+t \lambda(t)}{2+t \lambda(t)} =
\frac{1}{1 + \frac{\sqrt{t}}{g''(\sqrt{t})} g'(\sqrt{t})} =: \Theta(t)
\]
and remark that for $t \gg 1$ it holds
\[
c_1 t \ln\big(1+\sqrt{t}\big) \leq \frac{\sqrt{t}g'(\sqrt{t})}{g''(\sqrt{t})} \leq c_2 t \ln \big(1+\sqrt{t}\big)
\]
as well as $t \leq c_3 t \ln \big(1+\sqrt{t}\big)$, hence
\[
\Theta(t) \geq c_4 \frac{1}{t \ln\big(1+\sqrt{t}\big)} \geq c_5 \frac{1}{t \ln\big(1+t\big)}
\]
again for $t \gg 1$. Since
\[
\int_1^\infty \frac{\D t}{t\ln\big(1+t\big)} = \infty\, ,
\]
the claim \eqref{app 1} follows. \hspace*{\fill}\qed\\

\begin{remark}\label{app rem 1}
Formally the above model case of nearly linear growth should satisfy the $C^{3,\alpha}$-condition posed in \cite{NN:1959_1}. 
Since the integral occurring in \eqref{app 1} is not depending on the energy density for small values of $t$, we may easily adjust the example
to obtain a $C^{3,\alpha}$-density of nearly linear growth. For example we just let
$f_\eps(t) := \sqrt{\eps +t} \ln(1+\sqrt{\eps +t})$, $t\geq 0$, with some $\eps > 0$.
\end{remark}

\bibliography{bernsteinvariante}
\bibliographystyle{unsrt}

\vspace*{0.5cm}
\begin{tabular}{ll}
Michael Bildhauer&bibi@math.uni-sb.de\\
Martin Fuchs&fuchs@math.uni-sb.de\\[5ex]
Department of Mathematics&\\
Saarland University&\\
P.O.~Box 15 11 50&\\
66041 Saarbr\"ucken&\\ 
Germany&
\end{tabular}
\end{document}